\documentclass{aims} 
\usepackage{amsmath, amsfonts}
\usepackage{mathtools}
\usepackage{xcolor}
\usepackage{paralist}
\usepackage[misc]{ifsym}
\usepackage{epsfig} 
\usepackage{epstopdf} 
\usepackage[colorlinks=true]{hyperref}
\usepackage{listings}
\usepackage{courier}
\usepackage{layouts}
\usepackage{booktabs}
 
\definecolor{codegreen}{rgb}{0,0.4,0}
\definecolor{codegray}{rgb}{0.5,0.5,0.5}
\definecolor{codepurple}{rgb}{0.58,0,0.12}
\definecolor{backcolour}{rgb}{1.0,0.99,0.95}
\definecolor{deepblue}{rgb}{0.,0.1,0.8}
\definecolor{deepred}{rgb}{0.6,0,0}
\definecolor{deepgreen}{rgb}{0,0.5,0}
 
\lstdefinestyle{mystyle}{
    backgroundcolor=\color{backcolour},   
    commentstyle=\color{codegreen},
    keywordstyle=\color{deepblue},
    numberstyle=\tiny\color{codegray},
    stringstyle=\color{codepurple},
    breakatwhitespace=false,         
    breaklines=true,                 
    captionpos=b,                    
    keepspaces=true,                 
    showspaces=false,                
    showstringspaces=false,
    showtabs=false,                  
    tabsize=2,
    aboveskip=1em,
    belowskip=1em
}
 
\hypersetup{urlcolor=blue, citecolor=red}
\allowdisplaybreaks

\def\currentvolume{X}
 \def\currentissue{X}
  \def\currentyear{200X}
   \def\currentmonth{XX}
    \def\ppages{X-XX}
     \def\DOI{10.3934/jcd.xxxxx}

\newtheorem{theorem}{Theorem}[section]

\newtheorem{lemma}[theorem]{Lemma}

\theoremstyle{definition}

\newtheorem{remark}[theorem]{Remark}

\newcommand{\code}[1]{\lstinline$#1$}

\newcommand{\Cc}{\mathbb{C}}
\newcommand{\Gg}{\mathfrak{g}}

\title[ISOSYRK as Lie group methods]
{Isospectral symplectic Runge--Kutta schemes as Lie group methods} 

\author[Paolo Cifani, Klas Modin, Cecilia Pagliantini and Milo Viviani]{}


\thanks{The authors would like to thank all the participants of the mini-workshop  ``Geometric and Stochastic Methods for Fluid Models'', Kristineberg, Sweden, October 2024.}

\thanks{$^*$Corresponding author: Milo Viviani}

\begin{document}
\maketitle

\centerline{\scshape
Paolo Cifani$^{{\href{mailto:klas.modin@chalmers.se}{\textrm{\Letter}}}1}$
Klas Modin$^{{\href{mailto:klas.modin@chalmers.se}{\textrm{\Letter}}}2}$
Cecilia Pagliantini$^{{\href{mailto:klas.modin@chalmers.se}{\textrm{\Letter}}}3}$
and Milo Viviani$^{{\href{mailto:milo.viviani@sns.it}{\textrm{\Letter}}}*1}$
}

\medskip 

{\footnotesize
 \centerline{$^1$Scuola Normale Superiore Pisa, Italy}
}


{\footnotesize
 \centerline{$^2$Chalmers University of Technology and University of Gothenburg, Sweden}
}


{\footnotesize
 \centerline{$^3$University of Pisa, Italy}
} 


\bigskip




\begin{abstract}
We compare three approaches for structure preserving numerical integration of isospectral flows on quadratic Lie algebras.
Such flows originate from Hamiltonian dynamics on the cotangent bundle of the Lie group.
It is known, via discrete reduction theory, that symplectic Runge--Kutta methods applied to the cotangent bundle formulation induce isospectral symplectic Runge--Kutta (ISOSYRK) schemes on the Lie algebra.
Here, we show that the same symplectic Runge--Kutta method, but applied to the transport formulation of the flow on the Lie group, is equivalent to the corresponding ISOSYRK scheme.
We also give numerical results suggesting that the formulation on the Lie group is more efficient for schemes with two or more intermediate stages.
\end{abstract}


\section{Introduction}
Lie--Poisson systems are ubiquitous in mathematical physics. 
They arise from Hamiltonian mechanics when the phase space is the cotangent bundle $T^*G$ of a Lie group~$G$ and the symmetry of the Hamiltonian is the entire configuration space~$G$.
Due to the symmetry, the effective phase space then becomes the quotient $T^*G/G$, which via left or right translation is identified with the dual~$\mathfrak{g}^*$ of the Lie algebra~$\mathfrak{g}$ of $G$.
Examples include the rigid body, the spinning top, and the Heisenberg spin chain (see the monograph by Marsden and Ratiu~\cite{MaRa1999} for a comprehensive introduction and more examples).
The motivating example here is Zeitlin's model (\emph{cf.}~\cite{Ze1991,Ze2004,MoVi2020b}), which arises as a spatial Lie--Poisson discretization of the two-dimensional incompressible Euler equations.
However, the result we give applies more generally.

Indeed, we consider the case when the Lie group $G$ is quadratic, namely for a matrix $J\in \mathrm{GL}(n,\mathbb{C})$ its Lie algebra is given by
\begin{equation*}
    \mathfrak{g} = \{W \in \mathrm{GL}(n,\mathbb{C})\mid W^\dagger J + J W = 0 \}.
\end{equation*}
This setting includes, for example, the unitary, orthogonal, and symplectic groups. 
By identifying $\mathfrak{g}^*$ with $\mathfrak{g}$ via the Frobenius inner product, Lie--Poisson systems in this quadratic case take the form
\begin{equation}\label{eq:isospectral}
    \dot W = [\underbrace{\mathrm d H(W)^\dagger}_{B(W)}, W]
\end{equation}
where $H\colon \mathfrak{g}\to \mathbb{R}$ is the Hamiltonian and $[\cdot,\cdot]$ is the matrix commutator (see ~\cite{ModViv2020} for the full derivation of these equations from the abstract Lie--Poisson equations).
The equations~\eqref{eq:isospectral} are isospectral, i.e., the spectrum of $W$ is preserved by the flow.
From the point-of-view of Lie--Poisson dynamics, the isospectral property reflects that the flow restricts to \emph{coadjoint orbits} (see~\cite{MaRa1999} for details).
However, this feature is only the first aspect of the dynamics.
The second aspect is more delicate: each coadjoint orbit has a natural symplectic structure, the Kirillov--Kostant form, which is preserved by the flow.

For the numerical time integration of Lie--Poisson systems, it is usually desirable to preserve the geometric structure of the exact flow, especially in long-time simulations (\emph{cf.}~geometric numerical integration~\cite{HaLuWa2006}).
Indeed, the aim is to construct time-integration schemes whose discrete flow preserve the coadjoint orbits and are symplectic with respect to the Kirillov--Kostant form.
Whereas the former property is direct to verify, the second one is more subtle.
For a general Lie algebra $\mathfrak{g}$, a class of Lie--Poisson integrators is given by Bogfjellmo and Marthinsen~\cite{BoMa2016} based on Runge--Kutta--Munthe--Kaas and Crouch--Grossman type methods.
However, these methods rely on the exponential map, which makes them inefficient for high-dimensional Lie algebras (such as Zeitlin's model, see below).
In the case of quadratic Lie groups, an option is to use the isospectral symplectic Runge--Kutta (ISOSYRK) schemes, developed by Modin and Viviani~\cite{ModViv2020}.
These methods are free of exponential maps.
They are formulated on the level of the Lie algebra, i.e., on the level of the isospectral flow~\eqref{eq:isospectral}.
However, they are not the symplectic Runge--Kutta methods applied directly to the isospectral flow (which fail to preserve isospectrality).
Instead, they originate from a re-interpretation of the Runge--Kutta tableau in the context of Lie-Poisson reduction, which gives a new class of methods.
The relation to the corresponding standard symplectic Runge--Kutta method is obtained via un-reduction: in the un-reduced, canonical variables $(Q,P)\in T^*G$ the method is the standard symplectic Runge--Kutta method.

The purpose of this note is twofold.
As a main result, in Section~\ref{sec:isonum}, we prove that the ISOSYRK methods, when expressed directly on the Lie group $G$, are precisely given as the corresponding symplectic Runge--Kutta method applied to this Lie group reformulation (see equation~\eqref{eq:Q_dynamics_ql} below).
That the discrete flow obtained this way is isospectral has been known since long; it is a direct consequence of the fact that symplectic Runge--Kutta methods preserve quadratic invariants~\cite{Sa1988}.
However, since $G$ is not a symplectic manifold, it is surprising that the methods so obtained are precisely corresponding to the ISOSYRK methods and are therefore symplectic on the coadjoint orbits.

The second part, in Section~\ref{sec:numerics}, compares the efficiency of the formulations on $\mathfrak{g}$, $T^* G$, and $G$ in the case of Zeitlin's model (where $G = \mathrm{SU}(n)$ for a large $n$).
The conclusion is that when the number of intermediate stages $s$ exceeds 1, the most efficient formulation is on $G$, but for $s=1$ the most efficient formulation is on~$\mathfrak{g}$.

\section{Equivalent formulations of isospectral flows on quadratic Lie algebras}

In this section, we show three different equivalent formulations of isospectral flows on semisimple quadratic Lie algebras.
Consider a quadratic Lie group $G:=\lbrace Q\in GL(n,\Cc)\mid Q^\dagger JQ=J\rbrace$, for some $J\in GL(n,\Cc)$.
Its Lie algebra is defined as $\mathfrak{g}=\lbrace W\in \mathbb{M}(n,\Cc)\mid W^\dagger J+JW=0\rbrace$.
The group $G$ and its algebra $\Gg$ are called $J$-quadratic Lie group and $J$-quadratic Lie algebra, respectively.
Isospectral flows on $\Gg$ are defined as
\begin{equation}\label{eq:isospectral_ql}
\dot{W}=[B(W),W],    
\end{equation}
for some $B:\mathfrak{g}\rightarrow\mathfrak{n}(\mathfrak{g}):=\lbrace A\in \mathbb{M}(n,\Cc)\mid[A,\Gg]\subset\Gg\rbrace$.

In view of solving equation \eqref{eq:isospectral_ql} numerically, we give two equivalent formulations which are crucial in order to retain the invariants at a discrete level. 
In particular, we can lift equation \eqref{eq:isospectral_ql} to the cotangent bundle of $G$ in the following way.
Identify $\mathfrak{g}^*$ with $\mathfrak{g}$ via the Frobenius inner product, and consider the momentum map $\mu:T^*G\rightarrow\Gg^*$ given by the canonical left-action of $G$ on $T^*G$ and explicitly defined by $(Q,P)\mapsto Q^\dagger P$.
In the $(Q,P)$ variables satisfying $Q^\dagger P=W$, assuming that both $B$ and $B^\dagger$ are in $\Gg$, equation \eqref{eq:isospectral_ql} is equivalent to the system
\begin{equation}\label{eq:hamilton_ql}
    \begin{array}{ll}
        \dot{Q} &= QB(Q^\dagger P)^\dagger\\
        \dot{P} &= -PB(Q^\dagger P),
    \end{array}
\end{equation}
for any choice of $(Q_0,P_0)\in T^*G$ satifying $Q_0^{\dagger} P_0=W_0$.
In particular, for $Q_0=I_{n\times n}$ and $P_0=W_0$, it can be verified that the solution $(Q,P)$ to equation \eqref{eq:hamilton_ql} is such that $P=W_0 Q^{-\dagger}$. Hence, the equivalent solution $W$ to equation \eqref{eq:isospectral_ql} satisfies $W=Q^{\dagger}P=Q^\dagger W_0Q^{-\dagger}$. Then, using the fact that $Q\in G$, we have
\begin{equation}\label{eq:Q_dynamics_ql}
    \dot{Q} = QB(Q^\dagger W_0Q^{-\dagger})^\dagger=QB(Q^\dagger W_0JQJ^{-1})^\dagger.
\end{equation}
The assumption that both $B$ and $B^\dagger$ are in $\Gg$ depends only on $J$, as explained in the following lemma, proved in \cite{Viv2020}.
\begin{lemma}\label{theo:lemma} Let $\Gg$ be a Lie subalgebra of $\mathfrak{gl}(n,\Cc)$ such that there exists a matrix $J$ for which
\[W \in\Gg \iff W^\dagger J+JW=0.\]
Then $\Gg = \Gg^\dagger$ implies $J^2\in\mathfrak{c}(\Gg):=\lbrace A\mid[A,\Gg]=0\rbrace$. 
Moreover, if $J$ is invertible and $J^2\in\mathfrak{c}(\Gg)$, then $\Gg = \Gg^\dagger$.
\end{lemma}
In many cases, as for the irreducible representations of classical quadratic simple Lie algebras, $\mathfrak{c}(\Gg)$ consists only of a multiple of the identity. 
For this reason, in the next section, we consider only invertible $J$ such that $J^2$ is a non-null multiple of the identity matrix.
As shown in \cite[Prop. 6.28]{Kna1996}, the condition of being closed under adjunction is equivalent for the Lie algebra to be semisimple.
Hence, we have shown that the isospectral flow \eqref{eq:isospectral_ql} admits equivalent formulations \eqref{eq:hamilton_ql} and \eqref{eq:Q_dynamics_ql} when $\mathfrak{g}$ is a semisimple quadratic Lie algebra.
In these cases, without loss of generality, we can assume that $B:\Gg\rightarrow\Gg$.
In \cite{ModViv2020}, it is shown that when $B^\dagger$ is the gradient of some function $H$ on a semisimple Lie algebra $\Gg$, in the above identification of $\Gg^*$ and $\Gg$, the isospectral flow \eqref{eq:isospectral_ql} is a Lie--Poisson system on $\Gg^*$ and equations \eqref{eq:hamilton_ql} become the canonical Hamilton equations on $T^*G$.

\section{Isospectral numerical integration}\label{sec:isonum}

In this section, we show that the isospectral symplectic Runge--Kutta methods derived in \cite{ModViv2020} can be simply restricted to the integration of \eqref{eq:Q_dynamics_ql} with a symplectic Runge--Kutta method, reducing the number of intermediate stages from $s^2+s$ to the optimal $s$.
We recall that so far only in the case of the isospectral midpoint \cite{Viv2020} and of SDIRK \cite{DaSLes2022} the number of intermediate stages was already optimal.

Let $G$ be a $J$-quadratic Lie group, such that $J^2=\alpha I$ for $\alpha\neq 0$.
In \cite[Corollary 1]{ModViv2020}, it is shown that any numerical method that preserves $G$, i.e., quadratic invariants, applied to equation \eqref{eq:hamilton_ql} induces an isospectral method for \eqref{eq:isospectral_ql}.
In the case of $B=dH^\dagger$, for some function $H$, a symplectic Runge--Kutta method applied to equation~\eqref{eq:hamilton_ql} induces a Lie--Poisson integrator for \eqref{eq:isospectral_ql}, see \cite[Theorem 3]{ModViv2020}.

The ISOSYRK methods \cite{ModViv2020} form a class of numerical integrators for \eqref{eq:isospectral_ql}, which are equivalent to symplectic Runge--Kutta methods applied to equation \eqref{eq:Q_dynamics_ql}.
They are defined from a symplectic Butcher tableu $(a_{ij}, b_i)$ via the intermediate stages $X_i, K_{ij}$ by the equations
\begin{equation}\label{eq:isosyrk}
\left\{ \begin{aligned} 
X_i &= - h(W_k+\sum_{j=1}^s a_{ij} X_j)B(\widetilde{W}_i)  && \mbox{for } i=1,...,s\\
K_{ij} &= hB(\widetilde{W}_i)(\sum_{j'=1}^s (a_{ij'}X_{j'}+a_{jj'}K_{ij'}))    &&\mbox{for } i,j=1,...,s\\
\widetilde{W}_i&=W_k+\sum_{j=1}^s a_{ij} (X_j-J^{-1}X_j^\dagger J+K_{ij})    &&\mbox{for } i=1,...,s\\
W_{k+1} &= W_k + h\sum_{i=1}^sb_i[B(\widetilde{W}_i),\widetilde{W}_i].&& 
\end{aligned}\right. 
\end{equation}
The isospectrality of \eqref{eq:isosyrk} relies on the conservation of quadratic invariants by the underlying Runge--Kutta method, which is equivalent to its symplecticity.
Here we show that the numerical solution $\lbrace (Q_k,P_k)\rbrace_{k\geq 1}$, obtained via the symplectic Runge--Kutta method of Butcher tableu $(a_{ij}, b_i)$ applied to equation \eqref{eq:Q_dynamics_ql}, coincides to $\lbrace (\hat{Q}_k,P_0\hat{Q}_k^{-\dagger})\rbrace_{k\geq 1}$, where $\lbrace \hat{Q}_k\rbrace_{k\geq 1}$ is obtained via the same symplectic Runge--Kutta method applied to equation \eqref{eq:hamilton_ql}. Therefore, we have obtained a Lie--Poisson integrator for \eqref{eq:isospectral_ql} via the map $W_0\mapsto W_k:=\hat{Q}_k W_0 J\hat{Q}_k J^{-1}$.
Indeed, let us define the intermediate stages of a Runge--Kutta method for \eqref{eq:hamilton_ql} as
\begin{equation}\label{eq:intermediate_stages_ql} 
\begin{array}{cc}
    K^Q_i= & Q_0 + h\sum_j a_{ij}K^Q_jB((K^Q_j)^\dagger K^P_j)^\dagger \\
    K^P_i= & P_0 - h\sum_j a_{ij}K^P_jB((K^Q_j)^\dagger K^P_j).
\end{array}
\end{equation}
and for \eqref{eq:Q_dynamics_ql}
\begin{equation}\label{eq:intermediate_stages2_ql}
    \hat{K}^Q_i=  Q_0 + h\sum_j a_{ij}\hat{K}^Q_jB((\hat{K}^Q_j)^\dagger W_0J\hat{K}^Q_j J^{-1})^\dagger,
\end{equation}
for $i,j=1,\dots s$. 
We see that taking $K^P_i= W_0 JK^Q_iJ^{-1}$ and $P_0= W_0 JQ_0J^{-1}$, and replacing in the second equation of \eqref{eq:intermediate_stages_ql}, we get 
\begin{equation}
\begin{array}{ll}
    W_0JK^Q_iJ^{-1}&=  W_0JQ_0J^{-1} - h\sum_j a_{ij} W_0 JK^Q_jJ^{-1}B((K^Q_j)^\dagger W_0 JK^Q_jJ^{-1})\\
    &= W_0JQ_0J^{-1} + h\sum_j a_{ij} W_0 JK^Q_j(J^2)^{-1}B((K^Q_j)^\dagger W_0 JK^Q_jJ^{-1})^\dagger J
\end{array}
\end{equation}
where we have used the fact that $-J^{-1}B^\dagger J=B$. 
Furthermore, since $J^2=\alpha I$, $\alpha\neq 0$, then we conclude that $(K^Q_i,K^P_i)$ are indeed a solution of \eqref{eq:intermediate_stages_ql} if and only if $K^Q_i=\hat{K}^Q_i$ is a solution of \eqref{eq:intermediate_stages2_ql} and $K^P_i= W_0 J\hat{K}^Q_iJ^{-1}$.
With these calcuations we arrive at the main theoretical result of this paper.
\begin{theorem}
    Let $\Gg$ be a $J$-quadratic Lie algebra for which $J^2=\alpha I$ for some $\alpha\neq 0$.
    Further, let $h>0$ and consider an $s$-stage symplectic Runge--Kutta method with tableau $A=[a_{ij}]$ and $b=[b_i]$ for $i,j=1,...,s$. 
    Then, for initial data $W_0\in\Gg$ and any $k\in\mathbb{N}$, the following methods $W_{k-1}\mapsto W_k$ coincide:
\begin{enumerate}
    \item $W_{k}$ given by the ISOSYRK-method \eqref{eq:isosyrk} applied to equation \eqref{eq:isospectral_ql}; 
    \item $W_k\coloneqq Q_k^\dagger P_k$ where $(Q_k,P_k)$ are given by the $s$-stage symplectic Runge--Kutta method applied to equation \eqref{eq:hamilton_ql} with initial data $Q_0=I$ and $P_0=W_0$;
    \item $W_k\coloneqq Q_k W_0 JQ_kJ^{-1}$ where $Q_k$ is given by the $s$-stage symplectic Runge--Kutta method applied to equation \eqref{eq:Q_dynamics_ql} with initial data $Q_0=I$.
\end{enumerate}
\end{theorem}

\begin{remark}
    For partitioned symplectic Runge--Kutta methods one does not have a corresponding equivalence (unless, of course, they are genuine symplectic Runge--Kutta methods).
\end{remark}

\section{Numerical tests}\label{sec:numerics}
In the last section, we compare the performances of symplectic Runge--Kutta methods applied to equation \eqref{eq:hamilton_ql} and \eqref{eq:Q_dynamics_ql} and of isospectral symplectic Runge--Kutta methods applied to equation \eqref{eq:isospectral_ql}.
We consider the Euler--Zeitlin equations defined in \cite{Ze2004}, which is an isospectral flow on $\mathfrak{su}(n)$ arising as a spatial discretization of the 2D Euler equations on a 2-sphere:
\begin{equation}
    \dot{W} = \frac{1}{\hbar}[\Delta_n^{-1}W,W].
\end{equation}
Here, $W\in\mathfrak{su}(n)$, $\hbar=2/\sqrt{n(n-1)}$, and the discrete Laplacian is defined by
\[\Delta_n W = \sum_{i=1}^3[S_i,[S_i,W]],\]
for $S_i$ a basis of the $n$-dimensional irreducible unitary representation of the Lie algebra $\mathfrak{so}(3)$.
For details on Zeitlin's model and matrix hydrodynamics, we refer to the publications~\cite{MoVi2024,MoVi2026b}.

In Table~1, we show the run time of $5$ seconds of simulation of the Euler--Zeitlin equations run in MATLAB, obtained by specifying the three different schemes of the previous section for the Gauss methods with $s$ intermediate stages, for $s=1,2,3$ (hence of order $2,4,6$), for $n=17, 33, 65, 129$.
For each $n$, the first line represents the Gauss methods applied to equation \eqref{eq:hamilton_ql}, the second line to equation \eqref{eq:Q_dynamics_ql} and the third line the corresponding ISOSYRK method~\eqref{eq:isosyrk} applied to equation \eqref{eq:isospectral_ql}.
In the numerical simulations, we use the convention of \cite{MoVi2024} for the Euler--Zeitlin equations.
We set the initial value $W_0\in\mathfrak{su}(n)$, activating only the modes $T^N_{lm}$, for $l\leq 4$ with the same values for each $N$, and such that the spectral norm $\|W_0\|=1$.
We set the time step equal to $h=0.1$ and the tolerance of the fixed-point iteration, needed to solve the implicit equations at the intermediate stages, to $10^{-13}$.

In the case of the midpoint method ($s=1$), the most efficient scheme is the isospectral one.\footnote{The same conclusion applies to the isospectral symplectic DIRK defined in \cite{DaSLes2022}, as they are just compositions of isospectral midpoint methods.}
However, for $s\geq 2$ the Gauss method applied to equation \eqref{eq:Q_dynamics_ql} outperforms the other two, with the isospectral one becoming really poor for large $n$ and $s\geq 3$.

\begin{table}
\begin{center}
\begin{tabular}{c c c c} 
\toprule $n$ & $s=1$ & $s=2$  & $s=3$ \\ 
\midrule
   & \hspace{.175cm}  0.094731  & 0.258834  &	  0.438856  \\
17 &\hspace{.175cm}    0.060177  & \textbf{0.142243}  &	  \textbf{0.161755} \\
   &\hspace{.175cm}   \textbf{0.051316}  & 0.415429  &	  0.633625  \\
\midrule
 & \hspace{.175cm}  0.810970  & 1.555489  &	  2.336420  \\
33 & \hspace{.175cm}  0.392726  & \textbf{0.919517}  &	  \textbf{1.322197}  \\
 & \hspace{.175cm}  \textbf{0.251208}  &  1.486655  &	  3.573453  \\
\midrule
 & \hspace{.175cm}   4.328540  &  8.054763  &	  11.126779 \\
65 & \hspace{.175cm}   2.395906  & \textbf{4.280282}  &	  \textbf{6.188609}  \\
  & \hspace{.175cm}  \textbf{1.117873}  & 6.662383  &	  15.233301  \\
\midrule
 & \hspace{0cm}   36.188823  & 52.035348  &	  60.215178  \\
129 & \hspace{0cm}   19.480669  & \textbf{26.844684}  &	  \textbf{31.284791}  \\
 & \hspace{.175cm}   \textbf{9.838273}  & 52.616202  &	  108.424245  \\ 
\bottomrule
\end{tabular}
\end{center}
\vspace{2ex}
\caption{Run time of $5$ seconds of simulation of the Euler--Zeitlin equations run in MATLAB and obtained by Gauss methods with $s$ intermediate stages, for $s=1,2,3$ and $n=17, 33, 65, 129$. For each $n$, the first and second line represent the Gauss methods applied to equation \eqref{eq:hamilton_ql} and \eqref{eq:Q_dynamics_ql} respectively, and the third line the corresponding isospectral integrator applied to equation \eqref{eq:isospectral_ql}.}\label{tbl:1}
\end{table}

To conclude this section, we have given numerical evidence that for Lie--Poisson flows on quadratic Lie algebras in high dimension, solving equation \eqref{eq:Q_dynamics_ql} with a symplectic Runge--Kutta scheme is preferable for integration orders higher than two.
The isospectral midpoint method remains the best choice for order two.
New implementations of ISOSYRK methods, such as a parallelized version, might make them more competitive also for high-order integration.
In general, to which extent high-order time-integration schemes are needed in fluid simulations is still not well understood; due to the turbulent and chaotic behavior in the dynamics, local accuracy quickly diminishes regardless of the order.









\bibliographystyle{AIMS}
\bibliography{bibliography}

\medskip
Received xxxx 20xx; revised xxxx 20xx; early access xxxx 20xx.
\medskip

\end{document}